\documentclass[secthm,seceqn,amsthm,ussrhead,12pt]{amsart}
\usepackage{amsmath,latexsym,comment}
\usepackage[english]{babel}
\usepackage[psamsfonts]{amssymb}
\usepackage{times}
\usepackage{cite}
\usepackage{pdflscape} 
\usepackage{ulem}
\usepackage[mathcal]{euscript}
\usepackage{tikz}
\usepackage{hyperref}
\usepackage{cancel}
\usetikzlibrary{arrows}

\newtheorem{Th}{Theorem}
\newtheorem{Lem}[Th]{Lemma}

\newtheorem{corollary}[Th]{Corollary}
\newtheorem{Def}[Th]{Definition}

\newenvironment{Proof}[1][Proof.]{\begin{trivlist}
\item[\hskip \labelsep {\bfseries #1}]}{\flushright
$\Box$\end{trivlist}}

\begin{document}
	\sloppy

%\vspace*{5cm}

{\Large Degenerations of Filippov  algebras
\footnote{The work was supported by   
CNPq 453845/2018-5;
FAPESP 18/15712-0;
RFBR 18-31-20004;
by the President's "Program Support of Young Russian Scientists" (grant MK-2262.2019.1).}
}

\medskip

\medskip

\medskip

\medskip
\textbf{Ivan Kaygorodov$^{a}$ \& Yury Volkov$^{b}$}
\medskip

{\tiny
$^{a}$ Universidade Federal do ABC, Santo Andr\'e, Brazil.

$^{b}$ Saint Petersburg state university, Saint Petersburg, Russia.
\smallskip

    E-mail addresses:\smallskip

    Ivan Kaygorodov (kaygorodov.ivan@gmail.com),
    
    Yury Volkov (wolf86\_666@list.ru).

}

       \vspace{0.3cm}

%{\bf Abstract.} 

We consider the variety of Filippov ($n$-Lie) algebra structures on an $(n+1)$-dimensional vector space. 
The group $GL_n(K)$ acts on it, and we study the orbit closures with respect to the Zariski topology. 
This leads to the definition of Filippov algebra degenerations. 
We present some fundamental results on such degenerations, including trace invariants and necessary degeneration criteria. 
Finally, we classify all orbit closures in the variety of complex $(n+1)$-dimensional Filippov $n$-ary algebras.

\

       \vspace{0.3cm}

       \vspace{0.3cm}

\section{Introduction}

       \vspace{0.3cm}

The contraction of a Lie algebra is a limiting process, which at first has been studied in physics \cite{13,10}.
For example, classical mechanics is a limiting case of quantum mechanics as $\hbar \to 0,$
described by a contraction of the Heisenberg-Weyl Lie algebras to the Abelian Lie algebra of the same dimension.
The study of contractions of binary algebras has a very long history (see, for example, \cite{c1,c2,c3,nesterpop}).
One of the first attempts to study the contractions of $n$-ary algebras occurred in \cite{deaz} where the authors considered the variety of Filippov algebras .

A more general definition than the contraction is often used in mathematics, the so-called degenerations. 
Here, one considers the variety of $n$-dimensional algebra structures and the closures of $GL_n(K)$-orbits with respect to the Zariski topology.
There are a lot of papers on degenerations of algebras (see, for example, \cite{BB09, BC99,GRH2, GRH3, ikv17, kv18, GRH}).
One of important problems in this direction is the description of so-called rigid algebras. 
These algebras are of big interest because the closures of their orbits under the action of a generalized linear group form irreducible components of a variety under consideration
(with respect to the Zariski topology). 
For example, rigid algebras were classified in the varieties of
low dimensional associative, Jordan, Lie and Leibniz algebras.
There are also works in which the full information about degenerations was found for some variety of algebras.
Here one can mention the descriptions of degenerations of low dimensional 
associative, Lie, pre-Lie, Malcev, Leibniz (see \cite{BB09,kpv,kppv})
and all $2$-dimensional (see \cite{kv17}) algebras.
There is an obvious connection between the notions of degeneration and deformation. If the algebra $\mathbb{A}$ degenerates to the algebra $\mathbb{B}$, then $\mathbb{B}$ can be deformed to $\mathbb{A}$ (even via a so-called jump deformation).

We want to point out that one of the sources for the notion of Filippov
algebras was the Nambu mechanics proposed in \cite{nambu}. 
Many connections are currently known between problems in the theory of $n$-ary algebras and mathematical
physics, in particular, between Nambu mechanics and Chern--Simons theory (see \cite{deaz2,BB1,BB2} for
more details).

The purposes of the current paper are to initiate the study of degenerations of $n$-ary algebras,
to give some important definitions and facts in this direction
and to apply all of this to describe all degenerations in the variety of $(n+1)$-dimensional $n$-ary Filippov algebras over $\mathbb{C}.$

%\newpage

\section{Preliminaries}

All spaces in this paper are over $\mathbb{C}$, and we write simply $dim$, $Hom$ and $\otimes$ instead of $dim_{\mathbb{C}}$, $Hom_{\mathbb{C}}$ and $\otimes_{\mathbb{C}}$. An $n$-ary algebra $\mathbb{A}$ is a vector space with an $n$-ary  operation, i.e. with an $n$-linear map from $\underbrace{\mathbb{A}\times \ldots \times \mathbb{A}}_{n\mbox{\scriptsize { times}}}$ to $\mathbb{A}$. We will denote this map by $[\cdot,\dots,\cdot]$, i.e. the result of application of our operation to $(x_1,\dots,x_n)\in \mathbb{A}^n$ is $[x_1,\dots,x_n]\in\mathbb{A}$. Throughout the paper, we fix some integer $n\ge 2$ and write simply algebra instead of $n$-ary algebra.

Given a $k$-dimensional vector space $\mathbb{V}$, the set $Hom(\mathbb{V}^{\otimes n},\mathbb{V}) \cong (\mathbb{V}^*)^{\otimes n}  \otimes \mathbb{V}$ 
is a vector space of dimension $k^{n+1}$. 
This space has a structure of the affine variety $\mathbb{C}^{k^{n+1}}.$ 
Indeed, let us fix a basis $e_1,\dots,e_k$ of $\mathbb{V}$. 
Then any $\mu\in Hom(\mathbb{V}^{\otimes n},\mathbb{V})$ is determined by $k^{n+1}$ structure constants $c_{i_1, \ldots, i_n}^j\in\mathbb{C}$ such that
$\mu(e_{i_1} \otimes \ldots \otimes e_{i_n})=\sum\limits_{j=1}^k c_{i_1, \ldots, i_n}^je_j$. 
A subset of $Hom(\mathbb{V}^{\otimes n},\mathbb{V})$ is {\it Zariski-closed} if it can be defined by a set of polynomial equations in the variables $c_{i_1, \ldots, i_n}^j$ ($1\le i_1, \ldots, i_n,j \le k$).

Let $\mathbb{T}$ be a set of polynomial identities.
All algebra structures on $\mathbb{V}$ satisfying the polynomial identities from $\mathbb{T}$ form a Zariski-closed subset of the variety $Hom(\mathbb{V}^{\otimes n},\mathbb{V})$. 
We denote this subset by $\mathbb{L(T)}$.
The general linear group $GL(\mathbb{V})$ acts on $\mathbb{L(T)}$ by conjugations:
$$ (g * \mu )(x_1\otimes \ldots \otimes x_n) = g\mu(g^{-1}x_1\otimes \ldots \otimes g^{-1}x_n)$$ 
for $x_1, \ldots, x_n\in \mathbb{V}$, $\mu\in \mathbb{L(T)}\subset Hom(\mathbb{V}^{\otimes n},\mathbb{V})$ and $g\in GL(\mathbb{V})$.
Thus, $\mathbb{L(T)}$ is decomposed into $GL(\mathbb{V})$-orbits that correspond to isomorphism classes of algebras. 
Let $O(\mu)$ denote the orbit of $\mu\in\mathbb{L(T)}$ under the action of $GL(\mathbb{V})$ and $\overline{X}$ denote the Zariski closure of a set $X\subset \mathbb{L(T)}$.

Let $\mathbb{A}$ and $\mathbb{B}$ be two $k$-dimensional algebras satisfying identities from $\mathbb{T}$ and $\mu,\lambda \in \mathbb{L(T)}$ represent $\mathbb{A}$ and $\mathbb{B}$ respectively.
We say that $\mathbb{A}$ degenerates to $\mathbb{B}$ and write $\mathbb{A}\to \mathbb{B}$ if $\lambda\in\overline{O(\mu)}$.
Note that in this case we have $\overline{O(\lambda)}\subset\overline{O(\mu)}$. 
Hence, the definition of a degeneration does not depend on the choice of $\mu$ and $\lambda$. 
If $\mathbb{A}\not\cong \mathbb{B}$, then the assertion $\mathbb{A}\to \mathbb{B}$ is called a {\it proper degeneration}. We write $\mathbb{A}\not\to \mathbb{B}$ if $\lambda\not\in\overline{O(\mu)}$ and call such an assertion {\it non-degeneration}.

Let $\mathbb{A}$ be represented by $\mu\in\mathbb{L(T)}$. Then  $\mathbb{A}$ is  {\it rigid} in $\mathbb{L(T)}$ if $O(\mu)$ is an open subset of $\mathbb{L(T)}$.
Let us recall that a subset of a variety is called irreducible if it can not be represented as a union of two non-trivial closed subsets. A maximal irreducible closed subset of a variety is called an {\it irreducible component}.
In particular, $A$ is rigid in $\mathbb{L(T)}$ if and only if $\overline{O(\mu)}$ is an irreducible component of $\mathbb{L(T)}$. It is well known that any affine variety is a finite union of its irreducible components.

In the present work we use the methods applied to Lie algebras in \cite{BC99,GRH,GRH2}.
First of all, if $\mathbb{A}\to \mathbb{B}$ and $\mathbb{A}\not\cong \mathbb{B}$, then $dim\,Aut(\mathbb{A})<dim\,Aut(\mathbb{B})$, where $Aut(\mathbb{A})$ is the automorphism group of $\mathbb{A}$.
It follows from the fact that if the structure $\mu\in Hom(\mathbb{V}^{\otimes n},\mathbb{V})$ represents $\mathbb{A}$, then the dimension of $O(\mu)$ equals $dim\,GL(\mathbb{V})-dim\,Aut(\mathbb{A})=k^2-Aut(\mathbb{A})$.
We will compute all the automorphism groups of the algebras under consideration and will check the assertions $\mathbb{A}\to \mathbb{B}$ only for such $\mathbb{A}$ and $\mathbb{B}$ that $dim\,Aut(\mathbb{A})<dim\,Aut(\mathbb{B})$. Secondly, if $\mathbb{A}\to \mathbb{D}$ and $\mathbb{D}\to \mathbb{B}$ then $\mathbb{A}\to \mathbb{B}$. If there is no $\mathbb{D}$ such that $\mathbb{A}\to \mathbb{D}$ and $\mathbb{D}\to \mathbb{B}$ are proper degenerations, then the assertion $\mathbb{A}\to \mathbb{B}$ is called a {\it primary degeneration}. If $dim\,Aut(\mathbb{A})<dim\,Aut(\mathbb{B})$ and there are no $\mathbb{D}$ and $\mathbb{E}$ such that $\mathbb{D}\to \mathbb{A}$, $\mathbb{B}\to \mathbb{E}$, $\mathbb{D}\not\to \mathbb{E}$ and one of the assertions $\mathbb{D}\to \mathbb{A}$ and $\mathbb{B}\to \mathbb{E}$ is a proper degeneration,  then the assertion $\mathbb{A} \not\to \mathbb{B}$ is called a {\it primary non-degeneration}. It suffices to prove only primary degenerations and non-degenerations to describe degenerations in the variety under consideration. It is easy to see that any algebra degenerates to the algebra $\mathbb{C}^k$ with zero multiplication. From now on we use this fact without mentioning it.

To prove primary degenerations, we will construct families of matrices parameterized by $t$. 
Namely, let $\mathbb{A}$ and $\mathbb{B}$ be two algebras represented by the structures $\mu$ and $\lambda$ from $\mathbb{L(T)}$ respectively. 
Let $e_1,\dots, e_k$ be a basis of $\mathbb{V}$ and $c_{i_1,\ldots,i_n}^j$ ($1\le i_1,\ldots,i_n,j\le k$) be the structure constants of $\lambda$ in this basis. 
If there exist $a_i^j(t)\in\mathbb{C}$ ($1\le i,j\le k$, $t\in\mathbb{C}^*$) such that $E_i^t=\sum\limits_{j=1}^ka_i^j(t)e_j$ ($1\le i\le k$) form a basis of $\mathbb{V}$ for any $t\in\mathbb{C}^*$, and the structure constants of $\mu$ in the basis $E_1^t,\dots, E_k^t$ are such polynomials $c_{i_1,\ldots,i_n}^j(t)\in\mathbb{C}[t]$ that $c_{i_1,\ldots,i_n}^j(0)=c_{i_1,\ldots,i_n}^j$, then $\mathbb{A}\to \mathbb{B}$. 
In this case  $E_1^t,\dots, E_k^t$ is called a {\it parameterized basis} for $\mathbb{A}\to \mathbb{B}$.

%Note also the following fact. Let $\mathbb{B}(\alpha)$ be a series of algebras parameterized by $\alpha\in\mathbb{C}$ and $e_1,\dots,e_k$ be a basis of $\mathbb{V}$. Suppose also that, for any $\alpha\in\mathbb{C}$, the algebra $\mathbb{B}(\alpha)$ can be represented by a structure $\mu(\alpha)\in\mathbb{L(T)}$ having structure constants $c_{i_1,\ldots, i_n}^j(\alpha)\in\mathbb{C}$ in the basis $e_1,\dots,e_k$, where $c_{i_1,\ldots,i_n}^j(t)\in\mathbb{C}[t]$ for all $1\le i_1,\ldots, i_n,j\le k$. Let $\mathbb{A}$ be an algebra such that $\mathbb{A}\to \mathbb{B}(\alpha)$ for $\alpha\in\mathbb{C}\setminus S$, where $S$ is a finite subset of $\mathbb{C}$. Then $\mathbb{A}\to \mathbb{B}(\alpha)$ for all $\alpha\in\mathbb{C}$. This follows from the fact that $\mathbb{B}(\alpha)\in \overline{ \{O(\mathbb{B}(\beta)) \}_{\beta\in\mathbb{C}\setminus S}}$  for any $\alpha\in\mathbb{C}$. Thus, to prove that $\mathbb{A}\to \mathbb{B}(\alpha)$ for all $\alpha\in\mathbb{C}$ it is enough to construct degenerations that are valid for all but finitely many $\alpha$.

Let us describe the methods for proving primary non-degenerations. The main tool for this is the following lemma.

\begin{Lem}[\cite{GRH}]\label{main}
Let $\mathcal{B}$ be a Borel subgroup of $GL(\mathbb{V})$ and $\mathcal{R}\subset \mathbb{L(T)}$ be a $\mathcal{B}$-stable closed subset.
If $\mathbb{A} \to \mathbb{B}$ and $\mathbb{A}$ can be represented by $\mu\in\mathcal{R}$ then there is $\lambda\in \mathcal{R}$ that represents $\mathbb{B}$.
\end{Lem}

This lemma was proved in \cite{GRH} for binary algebras, but the proof does not use the algebra structure at all, and thus it is a general fact about orbit closures.
Since any Borel subgroup of $GL_k(\mathbb{C})$ is conjugated with the subgroup of lower triangular matrices, Lemma \ref{main} can be applied in the following way.  Suppose that $Q$ is a system of polynomial equations in variables $x_{i_1,\ldots, i_n}^j$ ($1\le i_1, \ldots, i_n,j\le k$).  
Suppose that if $x_{i_1, \ldots, i_n}^j=c_{i_1, \ldots, i_n}^j$ ($1\le i_1,\ldots, i_n,j\le k$) is a solution of $Q$, then $x_{i_1, \ldots, i_n}^j=\tilde c_{i_1, \ldots, i_n}^j$ ($1\le i_1,\ldots, i_n,j\le k$) is a solution of $Q$ too in the following cases:
\begin{enumerate}
    \item $\tilde c_{i_1, \ldots, i_n}^j=\frac{\alpha_{i_1} \ldots \alpha_{i_n}}{\alpha_j}c_{i_1, \ldots, i_n}^j$ for some $\alpha_m\in\mathbb{C}^*$ ($1\le m\le k$);

\item there are some numbers $1\le u<v\le k$ and some $\alpha\in\mathbb{C}$ such that  
    $$
    \tilde c_{i_1,\ldots,i_n}^j= 
\sum\limits_{U\subset\{1,\dots,n\}}\left(\prod\limits_{t\in U}\delta_{i_t,u}\right)\alpha^{|U|}\left( c_{s_v^U(i_1,\dots,i_n)}^j-\delta_{j,v}\alpha c_{s_v^U(i_1,\dots,i_n)}^u\right).
    $$
\end{enumerate}
 Here $s_v^U:\{1,\dots,k\}^n\rightarrow \{1,\dots,k\}^n$ sends the $n$-tuple $(i_1,\dots,i_n)$ to the $n$-tuple $(j_1,\dots,j_n)$, where $j_t=i_t$ for $t\not\in U$ and $j_t=v$ for $t\in U$.
We denote by $\delta$ the Kronecker delta function, i.e. we set $\delta_{i,j}:=0$ for $i\not=j$ and $\delta_{i,i}:=1$ for any $i$. As usually, $|U|$ denotes the number of elements of the set $U$.
Let $\mathbb{A}$ and $\mathbb{B}$ be two algebras and let $\mu,\lambda$ be the structures in $\mathbb{L(T)}$ representing $\mathbb{A}$ and $\mathbb{B}$ respectively.
Assume that there is a basis $f_1,\dots,f_k$ of $\mathbb{V}$ such that the structure constants of $\mu$ in this basis form a solution to all equations in $Q$.
If the structure constants of $\lambda$ in any basis do not form a solution to all equations in $Q$, then $\mathbb{A}\not\to \mathbb{B}$.

Sometimes it is more convenient to apply some corollaries of Lemma \ref{main} instead of its direct use. We will need some auxiliary notation to formulate other criteria for non-degenerations.
For $1 \le i_1< i_2< \dots < i_{n-t} \le n$, $f^1,\dots,f^{n-t}\in \mathbb{V}^*$ and $a_1,\dots,a_n\in \mathbb{V}$, we denote by $conv_{1,\dots,{n-t}}^{i_1,\dots,i_{n-t}}(f^1\otimes\dots\otimes f^{n-t}\otimes a_1\otimes \dots\otimes a_n)$ the convolution of $f^1\otimes\dots\otimes f^{n-t}\otimes a_1\otimes \dots\otimes a_n$ with respect to the pairs of indices $(1,{n-t}+i_1),\dots,({n-t},{n-t}+i_{n-t})$, i.e. the tensor
$$f^1(a_{i_1})\dots f^{n-t}(a_{i_{n-t}})a_1\otimes\dots\otimes a_{i_1-1}\otimes a_{i_1+1}\otimes\dots\otimes a_{i_{n-t}-1}\otimes a_{i_{n-t}+1}\otimes\dots\otimes a_n\in \mathbb{V}^{\otimes t}$$
that can be obtained from $a_1\otimes \dots\otimes a_n$ by omitting the components at places $i_1,\dots,i_{n-t}$ and multiplying by the coefficient $f^1(a_{i_1})\dots f^{n-t}(a_{i_{n-t}})\in \mathbb{C}$. One can extend $conv_{1,\dots,{n-t}}^{i_1,\dots,i_{n-t}}$ to the linear map $$conv_{1,\dots,{n-t}}^{i_1,\dots,i_{n-t}}:(\mathbb{V}^*)^{\otimes {n-t}}\otimes \mathbb{V}^{\otimes n}\rightarrow \mathbb{V}^{\otimes t}.$$ For $T\in \mathbb{V}^{\otimes n}$ and $I=\{1,\dots,n\}\setminus \{i_1, \dots, i_{n-t}\}$, we denote by $supp_I(T)$ the subspace of $\mathbb{V}^{\otimes t}$ generate by the elements $conv_{1,\dots,{n-t}}^{i_1,\dots,i_{n-t}}(f^1\otimes\dots\otimes f^{n-t}\otimes T)$ for all $f^1,\dots,f^{n-t}\in \mathbb{V}^*$.
For an element $\sigma$ of the symmetric group $S_n$ on $n$ elements, we define also $\sigma:V^{\otimes n}\rightarrow V^{\otimes n}$ as the unique linear map such that $\sigma(a_1\otimes\dots\otimes a_n)=a_{\sigma(1)}\otimes\dots\otimes a_{\sigma(n)}$.
\begin{enumerate}
\item For $I\subset \{1,\dots,n\}$, $|I|=t$ we define the $I$-annihilator of $\mu\in\mathbb{L}(\mathbb{T})$ as the set 
$$Ann_I(\mu)=\{ X \in \mathbb{V}^{ \otimes t} \mid \mu(T)=0\mbox{ whenever $supp_I(T)$ is generated by $X$ over $\mathbb{C}$}\}.$$
It is not difficult to show that $Ann_I(\mu)$ is a linear subspace of $\mathbb{V}^{ \otimes t}$. This definition can be transferred to the algebra settings to define the $I$-annihilator of $\mathbb{A}$ as the corresponding subspace of $\mathbb{A}^{\otimes t}$.
\item $Ann(\mathbb{A})= \cap_{i=1}^n Ann_{\{i\} } (\mathbb{A})\subset\mathbb{A}$ is the  annihilator of $\mathbb{A}.$
%\item $msub_0(\mathbb{A})$ is a trivial subalgebra of $\mathbb{A}$ of the maximal dimension (we fix one for each algebra $\mathbb{A}$).
\item  Let us pick some $1\le t \le n$. By $S_{t,n-t}$ we denote the subgroup of the symmetric group $S_n$ on $n$ elements formed by such $\sigma$ that $\sigma(i)<\sigma(j)$ if either $1\le i<j\le t$ or $t+1\le i<j\le n$.
We define  the $t$-center $Z_t(\mu)$ of $\mu\in\mathbb{L}(\mathbb{T})$  as the set
$$Z_t(\mathbb{A})=\{  X \in \mathbb{V}^{ \otimes t} \mid \mu\sigma(X\otimes Y)=\mu(X\otimes Y)\mbox{ for any }Y\in\mathbb{V}^{ \otimes (n-t)}\mbox{ and any }
\sigma\in S_{t,n-t}  \}.$$ 
It is clear that $Z_t(\mu)$ is a linear subspace of $\mathbb{V}^{ \otimes t}$. This definition can be transferred to the algebra settings to define the $t$-center of $\mathbb{A}$ as the corresponding subspace of $\mathbb{A}^{\otimes t}$. In particular $Z(A)=Z_1(A)$ is the {\it center} of the $n$-ary algebra $A$.
\item By $\mathbb{A}^2$ we denote the subspace of $\mathbb{A}$ generated by $[a_1, \ldots, a_n]$ for all $a_1,\dots,a_n\in\mathbb{A}$.
%  $\{ a \in A \mid xa =ax, \mbox{ for all } x\in A \}.$
\end{enumerate}

Given linear spaces $\mathbb{U}$ and $\mathbb{W}$, we write simply $\mathbb{U}>\mathbb{W}$ instead of $dim\,\mathbb{U}>dim\,\mathbb{W}$. 
%We use also the notation $U\circ W=UW+WU$.
Then $\mathbb{A}\not\to \mathbb{B}$ in the following cases:
\begin{enumerate}
\item $Ann_I(\mathbb{A})>Ann_I(\mathbb{B})$ for some $I\subset \{1,\dots,n\}$;
\item $Ann(\mathbb{A})>Ann(\mathbb{B})$;
%\item $msub_0(\mathbb{A})>msub_0(\mathbb{B})$;;
\item $Z_t(\mathbb{A})>Z_t(\mathbb{B})$ for some $1\le t\le n$;
\item $\mathbb{A}^{2}<\mathbb{B}^{2}$.
\end{enumerate}
The second and the last criteria follow from Lemma \ref{main}. To show the first and the third criteria, it is enough to note that each of the conditions $X\in Ann_I(\mu)$ and $X\in Z_t(\mu)$ for $X\in\mathbb{V}^{\otimes t}$ is equivalent to a system of linear equations whose coefficients are polynomials in the structure constants of $\mu$. The fact that such a system has not less than $d$ linearly independent solutions for a given $d\ge 0$ is equivalent to the fact that its rank is less or equal to $(\dim\,\mathbb{V})^t-d$. Such a condition determines a closed subset $S_d$ of $\mathbb{L(T)}$ stable under the action of $GL(\mathbb{V})$. Since $\mathbb{A}$ can be represented by a structure from $S_d$ with $d=dim\,Ann_I(\mathbb{A})$ or $d=dim\,Z_t(\mathbb{A})$, depending on what we are proving, and $\mathbb{B}$ cannot be represented by such a structure, we have $\mathbb{A}\not\to \mathbb{B}$.

Following \cite{deaz}, we call an element of the form $X = (X_1,..., X_{n-1})\in \mathbb{A}^{n-1}$ a {\it fundamental object} for the $n$-ary algebra $\mathbb{A}$.
For a fundamental object $X$, one can define the right multiplication map (which in the case of Filippov and, more generally, right $n$-Leibniz algebras occurs to be an inner derivation) by the equality
$$R_XZ=Z\cdot X   := [Z, X_1, \ldots, X_{n-1}]\mbox{ for } Z \in \mathbb{A}.$$
The notion of a fundamental object allows to generalize the trace invariants that were introduced for binary algebras in \cite{BC99} to the case of $n$-ary algebras. 
Namely, if for $i,j \in \mathbb{N}$ there exists a unique constant $c_{i,j}(\mathbb{A})$ such that for any two fundamental objects  $X$ and $Y$ one has
\begin{equation}
\label{inv}{tr (R_X^i) tr  (R_Y^j)}=c_{i,j}(\mathbb{A}){tr (R_X^i \circ R_Y^j)},
\end{equation}
then $c_{i,j}(\mathbb{A})$ is called an {\it $(i,j)$-invariant} of $\mathbb{A}$. We also say that $\mathbb{A}$ has $(i,j)$-invariant $c_{i,j}(\mathbb{A})=\infty$ if $tr(R_X^i \circ R_Y^j)=0$ for any fundamental objects $X,Y$ and, at the same time, there are some $X,Y$ such that $tr (R_X^i) tr  (R_Y^j)\not=0$. Defined in this way $(i,j)$-invariants $c_{i,j}(\mathbb{A})\in\mathbb{C}\cup\{\infty\}$ satisfy the following property.

\begin{Lem}
Suppose that $\mathbb{A} \to \mathbb{B}$ and both $\mathbb{A}$ and $\mathbb{B}$ have $(i,j)$-invariants. Then $c_{i,j}(\mathbb{A})=c_{i,j}(\mathbb{B}).$
\end{Lem}

\begin{Proof}
The required assertion follows from the fact that the validness of the equality \eqref{inv} for all pairs $X,Y$ of fundamental objects can be expressed in terms of polynomial equations on the structure constants of $\mathbb{A}$.
\end{Proof}

\begin{Def}
For $(\alpha)=(\alpha_0, \alpha_1, \ldots, \alpha_n) \in \mathbb{C}^{n+1}$ let us define $Der_{(\alpha)}(\mathbb{A})$ as the space of all $D\in End(\mathbb{A})$
such that
$$\alpha_0 D[x_1, \ldots, x_n] = \sum\limits_{i=1}^n \alpha_i [x_1, \ldots, D(x_i), \ldots, x_n],$$
for all $x_1, \ldots, x_n \in \mathbb{A}.$
The elements of $Der_{(\alpha)} (\mathbb{A})$ are called $(\alpha)$-derivations.
\end{Def}

Note that if $(\alpha)=(1,1,\dots,1)$, then $(\alpha)$-derivation is the same thing as a usual derivation of an $n$-ary algebra. If $Der(\mathbb{A})$ denotes the set of derivation of the algebra $\mathbb{A}$, then $dim\,Der(\mathbb{A})=dim\,Aut(\mathbb{A})$.
Thus, one has $Der(\mathbb{A})<Der(\mathbb{B})$ in the case of a proper degeneration $\mathbb{A} \to \mathbb{B}$.

\begin{Lem}
If $\mathbb{A} \to \mathbb{B},$
then $Der_{(\alpha)}(\mathbb{A}) \leq Der_{(\alpha)}(\mathbb{B})$ for all $(\alpha) \in \mathbb{C}^{n+1}.$
\end{Lem}

\begin{Proof}
The proof is similar to the proof of \cite[Lemma 3.9]{BB09}.
\end{Proof}

In the cases where all of these criteria can not be applied to prove $\mathbb{A}\not\to \mathbb{B}$, we will define $\mathcal{R}$ by a set of polynomial equations and will give a basis of $\mathbb{V}$, in which the structure constants of $\mu$ representing $\mathbb{A}$ give a solution to all these equations. In this paper this will occur one time.
%We will omit everywhere the verification of the fact that $\mathcal{R}$ is stable under the action of the subgroup of lower triangular matrices and of the fact that $\lambda\not\in\mathcal{R}$ for any choice of a basis of $\mathbb{V}$. These verifications can be done by direct calculations.

If the number of orbits under the action of $GL(\mathbb{V})$ on  $\mathbb{L(T)}$ is finite, then the graph of primary degenerations gives the whole picture of the structure of the variety under consideration.
In particular,  irreducible components correspond in this case bijectively to rigid algebras and an algebra is rigid if and only if there is no proper degeneration to it.
If  $\mathbb{L(T)}$ contains infinitely many non-isomorphic algebras, we have to make some additional work to describe rigid algebras and irreducible components. For example, if the algebra $\mathbb{A}$ can be presented by a structure $\mu(0)$ contained in the family $\mu(x)\in\mathbb{L(T)}$ whose structure constants are polynomials in $x\in\mathbb{C}$ such that $\mu(x)\not\cong\mu(y)$ for $x\not=y$, then $\mathbb{A}$ is automatically not rigid while it is possible that there is no proper degeneration to it.

Let $\mathbb{A}(*):=\{\mathbb{A}(x)\}_{x\in X}$ be a set of algebras, and let $\mathbb{B}$ be another algebra. Suppose that, for $x\in X$, $\mathbb{A}(x)$ is represented by the structure $\mu(x)\in\mathbb{L(T)}$ and $\mathbb{B}\in\mathbb{L(T)}$ is represented by the structure $\lambda$. Then $\mathbb{A}(*)\to \mathbb{B}$ means $\lambda\in\overline{\{O(\mu(x))\}_{x\in X}}$, and $\mathbb{A}(*)\not\to \mathbb{B}$ means $\lambda\not\in\overline{\{O(\mu(x))\}_{x\in X}}$.

Let $\mathbb{A}(*)$, $\mathbb{B}$, $\mu(x)$ ($x\in X$) and $\lambda$ be as above. To prove $\mathbb{A}(*)\to \mathbb{B}$ it is enough to construct a family of pairs $(f(t), g(t))$ parameterized by $t\in\mathbb{C}^*$, where $f(t)\in X$ and $g(t)\in GL(\mathbb{V})$. Namely, let $e_1,\dots, e_k$ be a basis of $\mathbb{V}$ and $c_{i_1,\ldots, i_n}^j$ ($1\le i_1,\ldots, i_n,j \le k$) be the structure constants of $\lambda$ in this basis. 
If we construct $a_i^j:\mathbb{C}^*\to \mathbb{C}$ ($1\le i,j\le k$) and $f: \mathbb{C}^* \to X$ such that $E_i^t=\sum\limits_{j=1}^k a_i^j(t)e_j$ ($1\le i\le k$) form a basis of $\mathbb{V}$ for any  $t\in\mathbb{C}^*$, and the structure constants of $\mu_{f(t)}$ in the basis $E_1^t,\dots, E_k^t$ are such polynomials $c_{i_1, \ldots, i_n}^j(t)\in\mathbb{C}[t]$ that $c_{i_1, \ldots, i_n}^j(0)=c_{i_1, \ldots, i_n}^j$, 
then $\mathbb{A}(*)\to \mathbb{B}$. In this case  $E_1^t,\dots, E_k^t$ and $f(t)$ are called a parameterized basis and a {\it parameterized index} for $\mathbb{A}(*)\to \mathbb{B}$ respectively.

We now explain how to prove that $\mathbb{A}(*)\not\to \mathbb{B}$. Note first that if $dim\,Aut(\mathbb{A}(x))>dim\,Aut(\mathbb{B})$ for all $x\in X$, then $\mathbb{A}(*)\not\to \mathbb{B}$.
On the other hand, it is possible that $\mathbb{A}(*)\to \mathbb{B}$ and $dim\,Aut(\mathbb{A}(x))\ge dim\,Aut(\mathbb{B})$ for all $x\in X$.
One can use also the following generalization of Lemma \ref{main} that has the same proof.

\begin{Lem}\label{gmain}
Let $\mathcal{B}$ be a Borel subgroup of $GL(\mathbb{V})$ and $\mathcal{R}\subset \mathbb{L(T)}$ be a $\mathcal{B}$-stable closed subset.
If $\mathbb{A}(*) \to \mathbb{B}$ and for any $x\in X$ the algebra $\mathbb{A}(x)$ can be represented by a structure $\mu(x)\in\mathcal{R}$, then there is $\lambda\in \mathcal{R}$ representing $\mathbb{B}$.
\end{Lem}

In the same manner as before, one can show that $\mathbb{A}(*)\not\to \mathbb{B}$ in the following cases:
\begin{enumerate}
\item $Ann_I\big(\mathbb{A}(x)\big)>Ann_I(\mathbb{B})$ for some $I\subset \{1,\dots,n\}$ and all $x\in X$;
\item $Ann\big(\mathbb{A}(x)\big)>Ann(\mathbb{B})$ for  all $x\in X$;
%\item $msub_0(\mathbb{A})>msub_0(\mathbb{B})$;;
\item $Z_t\big(\mathbb{A}(x)\big)>Z_t(\mathbb{B})$ for some $1\le t\le n$ and all $x\in X$;
\item $\mathbb{A}(x)^{2}<\mathbb{B}^{2}$  for  all $x\in X$;
\item $\mathbb{A}(x)$ and $\mathbb{B}$ have $(i,j)$-invariants for all $x\in X$, $c_{i,j}\big(\mathbb{A}(x)\big)=c_{i,j}$ does bot depend on $x$ and $c_{i,j}\not=c_{i,j}(\mathbb{B})$;
\item $Der_{(\alpha)}\big(\mathbb{A}(x)\big) > Der_{(\alpha)}(\mathbb{B})$ for some $(\alpha) \in \mathbb{C}^{n+1}$ and all $x\in X$.
\end{enumerate}

%\newpage    
\section{The variety of Filippov algebras}

Filippov ($n$-Lie) algebras appeared in \cite{fil85} as a generalization of Lie algebras.
An algebra $\mathbb{A}$ with an $n$-ary multiplication $[ \cdot, \ldots, \cdot]:\mathbb{A}^n\rightarrow \mathbb{A}$ 
is called a {\it Filippov algebra} if the equalities
$$[x_1, \ldots, x_n] = (-1)^{\sigma}[x_{\sigma(1)}, \ldots, x_{\sigma(n)}]\,\, \forall\sigma \in S_n;$$
$$ [[x_1, \ldots, x_n], y_2, \ldots, y_n]= \sum\limits_{i=1}^n [x_1, \ldots, x_{i-1}, [x_i, y_2, \ldots, y_n], x_{i+1}, \ldots, x_n]$$
hold for all $x_1,\dots,x_n,y_2,\dots,y_n\in \mathbb{A}$. Here and later $(-1)^{\sigma}$ denotes the sign of the permutation $\sigma$.
We will denote by $\mathfrak{Fil}^n_k$ the variety of $k$-dimensional complex $n$-ary Filippov algebras. Note that $\mathfrak{Fil}^2_k$ is exactly the variety of $k$-dimensional Lie algebras.

Some properties of the variety of Filippov algebras are similar to the properties of the variety of Lie algebras.
One can define the solvable ideal of a Filippov algebra, simple and semisimple Filippov algebras, etc., see \cite{kas}.
Some properties of nilpotent Filippov algebras were studied in \cite{nil1,nil2,nil3}.
Two cohomological properties of semisimple Lie algebras also hold in the Filippov  algebras case. Namely, 
semisimple Filippov algebras are rigid (i.e. cannot be deformed in the Gerstenhaber sense) and do not admit non-trivial central extensions.
This result proved in \cite{deaz09} is an analogue of Whitehead's Lemma for Filippov algebras.
On the other hand, there are some properties of Lie algebras that Filippov $n$-ary algebras do not have for $n>2.$
For example, it is known that every finite-dimensional Lie algebra with an invertible derivation over a field of zero characteristic is nilpotent,
but it is not true for Filippov algebras \cite{wil}.

It is easy to see that up to isomorphism there is only one non-trivial anticommutative $n$-ary $n$-dimensional algebra.
All $(n+1)$-  and $(n+2)$-dimensional Filippov algebras over an algebraically closed field of characteristic zero were classified  \cite{bai11}.
Note that $(n+1)$-dimensional $n$-ary Filippov algebras play a crucial role in the classification of simple and semisimple Filippov algebras over an algebraically closed field of characteristic zero, because
in the case $n>2$ the $(n+1)$-dimensional algebra  $D_{n+1}$ from \cite{kac10} is the unique simple $n$-ary Filippov algebra. 
It is known that any finite-dimensional semisimple Filippov algebra over an algebraically closed field of characteristic zero is a direct sum of simple algebras. 
On the other hand, \cite{poj} gives many examples of simple Filippov algebras besides the algebra $D_{n+1}$ in the case of a field of positive characteristic.
Note also that there are some connections between Filippov algebras and Nambu-Poisson algebras, see \cite{kac16}.

In the table below we list up to isomorphism all algebras constituting the variety $\mathfrak{Fil}^n_{n+1}$ except the algebra $\mathbb{C}^{n+1}$ with zero multiplication.
This list is based on the classification presented in \cite[Lemma 3.1]{bai11}. In multiplication tables we give only values of nonzero products of the form $[e_{i_1},\dots,e_{i_n}]$ with $i_1<\dots<i_n$. All other products of basic elements are zero or can be obtained from the given ones via the equality $[e_{i_{\sigma(1)}},\dots,e_{i_{\sigma(n)}}]=(-1)^{\sigma}[e_{i_1},\dots,e_{i_n}]$ that holds for any $\sigma\in S_n$. The dimensions of automorphism groups are calculated in the next section.

\

\begin{center}
\begin{tabular}{|l|l|l|}
\hline

%$A_1$ & $[e_1, \ldots, \widehat{e}_i, \ldots, e_{n+1}]=0$ & \\
%\hline

$\mathbb{A}$ & {Multiplication table} & dim  Aut ($\mathbb{A}$) \\
\hline 
$B$ &  $[e_2, \ldots, e_{n+1}]=e_1 $  &$n^2+n$ \\
        
\hline

%$B_2$ &  $[e_1,  \ldots, e_{n}]=e_1$ & & $n^2$ \\
%\hline

$C_1$ &  $[e_2, \ldots,  e_{n+1}]=e_1,$ $[e_1, e_3,  \ldots, e_{n+1}]=e_2$ & $n^2$ \\
\hline

$C_2(\alpha)$ &  
$[e_2, \ldots,  e_{n+1}]=\alpha e_1+e_2,$ $[e_1, e_3,  \ldots, e_{n+1}]=e_2$ & $n^2$ \\
\hline

$C_3$ &  $[e_1, e_3, \ldots,  e_{n+1}]=e_1,$ $[e_2, e_3,  \ldots, e_{n+1}]=e_2$ & $n^2+2$ \\
\hline

$D_r, 3 \leq r \leq n+1$ &  $[e_1, \ldots, e_{i-1},e_{i+1}, \ldots, e_{n+1}]=e_i,$ $1 \leq i \leq r$ & $(n+1-r)(n+1)+\frac{r(r-1)}{2} $\\
\hline
\end{tabular}
\end{center}

\section{Automorphisms of $(n+1)$-dimensional $n$-ary Filippov algebras}

Let us recall the definition of an automorphism of an $n$-ary algebra.

\begin{Def}
Let $\mathbb{A}$ be an $n$-ary algebra. The linear map $\phi:\mathbb{A}\rightarrow \mathbb{A}$ is called an automorphisms of the algebra $\mathbb{A}$ if
$$\phi[x_1, \ldots, x_n]=[\phi(x_1), \ldots, \phi(x_n)]$$
for all $ x_1,\ldots, x_n \in \mathbb{A}.$
\end{Def}

Let us describe now the automorphisms of $(n+1)$-dimensional $n$-ary Filippov algebras. We will describe them by their matrices in the basis $e_1,\dots,e_{n+1}$.
Let $c_{i_1,\dots,i_n}^j$ ($1\le i_1,\dots,i_n,j\le n+1$) be the structure constants of the structure $\mu$ in the basis $e_1,\dots,e_{n+1}$. Then one can form a matrix $R$ of dimension $(n+1)\times(n+1)$ whose $(i,j)$-entry is $(-1)^{i-1}c_{1,\dots,i-1,i+1,\dots,n+1}^j$. The discussion before \cite[Theorem 2]{fil85} shows that the invertilble linear map having matrix $S$ in the basis $e_1,\dots,e_{n+1}$ determines an automorphism of  $\mu$ if and only if $\frac{1}{det(S)}SRS^T=R$.

%Let $\phi$ be an automorphism of an $n$-ary algebra $\mathbb{A}$ with the basis $\{ e_1, \ldots, e_m \},$ then 
%$\phi(e_i)=\sum \phi^i_j e_j$ and 
%$$\mathbb{A}_i=\sum\limits_{{\tiny \begin{array}{c} \sigma \in S_{n+1}, \\ \sigma(1)=1, \ldots, \sigma(i)=i \end{array} }} (-1)^{\sigma} \phi^{i+1}_{\sigma(i+1)} \ldots \phi^{n+1}_{\sigma(n+1)}.$$
%Below we are computing automorphisms for all algebras from $\mathfrak{Fil}^n_{n+1}.$

\begin{enumerate}
    \item Any automorphism of $B$ has the form
$$
\phi=\begin{pmatrix}
det(U)&a_2&\dots&a_{n+1}\\
0&&&\\
\vdots&&U&\\
0&&&
\end{pmatrix},
$$
where $U$ is an arbitrary invertible matrix of dimension $n\times n$ and $a_2,\dots,a_{n+1}$ are arbitrary complex numbers.
It is clear also that any linear map of this form is an automorphism of $B$.
In particular, one has $dim\,Aut(B) = n^2+n.$

\item Any automorphism of $C_1$ has the form
$$
\phi=\begin{pmatrix}
a&b\,det(U)&a_{1,3}&\dots&a_{1,n+1}\\
b&a\,det(U)&a_{2,3}&\dots&a_{2,n+1}\\
0&0&&&\\
\vdots&\vdots&&U&\\
0&0&&&
\end{pmatrix},
$$
where $U$ is a matrix of dimension $(n-1)\times (n-1)$ such that $det(U)=\pm 1$, $a,b$ are complex numbers such that $a^2\not=b^2$ and  $a_{1,3},\dots,a_{1,n+1},a_{2,3},\dots,a_{2,n+1}$ are arbitrary complex numbers.
It is clear also that any linear map of this form is an automorphism of $C_1$.
In particular, one has
$$dim\,Aut(C_1) = (n-1)^2-1+2+2(n-1)=n^2.$$

\item Any automorphism of $C_2(\alpha)$ has the form
$$
\phi=\begin{pmatrix}
a&\alpha b&a_{1,3}&\dots&a_{1,n+1}\\
b&a+b&a_{2,3}&\dots&a_{2,n+1}\\
0&0&&&\\
\vdots&\vdots&&U&\\
0&0&&&
\end{pmatrix},
$$
where $U$ is a matrix of dimension $(n-1)\times (n-1)$  such that $det(U)=1$, $a,b$ are complex numbers such that $a(a+b)\not=\alpha b^2$ and  $a_{1,3},\dots,a_{1,n+1},a_{2,3},\dots,a_{2,n+1}$ are arbitrary complex numbers.
It is clear also that any linear map of this form is an automorphism of $C_2(\alpha)$.
In particular, one has
$$dim\,Aut\big(C_2(\alpha)\big) = (n-1)^2-1+2+2(n-1)=n^2.$$

\item Any automorphism of $C_3$ has the form
$
\phi=\begin{pmatrix}
U&V\\
0&W
\end{pmatrix},
$
where $U$  is an arbitrary invertible matrix of dimension $2\times 2$, $W$ is a matrix of dimension $(n-1)\times (n-1)$ such that $det(W)=1$ and $V$ is an arbitrary matrix of dimension $2\times (n-1)$. It is clear also that any linear map of this form is an automorphism of $C_3$.
In particular, one has
$$dim\,Aut(C_3) = 4+(n-1)^2-1+2(n-1)=n^2+2.$$

    \item Any automorphism of $D_r$ has the form
$
\phi=\begin{pmatrix}
U&V\\
0&W
\end{pmatrix},
$
where $U$  is an invertible matrix of dimension $r\times r$, $W$ is an invertible matrix of dimension $(n+1-r)\times (n+1-r)$ and $V$ is an arbitrary matrix of dimension $r\times (n+1-r)$.
Due to the remark above, the map $\phi$ of such form gives an automorphism of $D_r$ if and only if
$$
\frac{1}{det(U)det(W)}U S_rU^T=S_r,
$$
where $S_r=diag(1,-1,\dots,(-1)^{r-1})$ is the diagonal matrix of dimension $r\times r$ whose $(i,i)$-entry is equal to $(-1)^{i-1}$. Hence, $\phi$ is an automorphism if and only if $det(U)^{r-2}det(W)^{r}=1$ and the matrix $Y_r^{-1}\frac{U}{a}Y_r$ is orthogonal, where $Y_r=diag(1,\mathrm{i},\dots,\mathrm{i}^{(r-1)^2})$ is a diagonal matrix such that $Y_r^2=S_r$ and $a$ is a complex number such that $a^2=\det(U)\det(W)$.

Thus, the automorphisms of $D_r$ are linear maps of the form $
\phi=\begin{pmatrix}
aY_rUY_r^{-1}&V\\
0&W
\end{pmatrix},
$
where $U$ is an orthogonal matrix of dimension $r\times r$, $W$ is an arbitrary invertible matrix of dimension $(n+1-r)\times (n+1-r)$, $V$ is an arbitrary matrix of dimension $r\times (n+1-r)$ and $a$ is a complex number such that
$a^{r-2}=\frac{1}{det(U)det(W)}$.
Since the orthogonal group of $r\times r$ matrices has dimension $\frac{r(r-1)}{2}$, one has
$$dim\,Aut(D_r) = \frac{r(r-1)}{2}+(n+1-r)^2+r(n+1-r)=(n+1-r)(n+1)+\frac{r(r-1)}{2}.$$
\end{enumerate}

%\newpage 
\section{Degenerations of $(n+1)$-dimensional $n$-ary Filippov algebras}

%Two cohomological properties of semisimple Lie algebras also hold for Filippov algebras, namely, that semisimple Filippov algebras do not admit non-trivial central extensions and that they are rigid i.e., cannot be deformed in Gerstenhaber sense. 
At this moment there are a few papers devoted to deformations and degenerations of $n$-ary algebras.
Nevertheless, there are some paper on this topic, for example, \cite{mak07} and \cite{mak16} can be mentioned among them.
The problem of description of degenerations in the varieties of $(n+1)-$ and $(n+2)-$dimensional $n$-ary Filippov algebras was posted in the second of the mentioned papers.
In this paper we describe all degenerations, rigid algebras and irreducible components in the variety $\mathfrak{Fil}^n_{n+1}$, and hence solve one of the problems posted in \cite{mak16}.

\begin{Th}\label{theorem}
The graph of primary degenerations for $\mathfrak{Fil}^n_{n+1}$ has the following form: 
\end{Th}
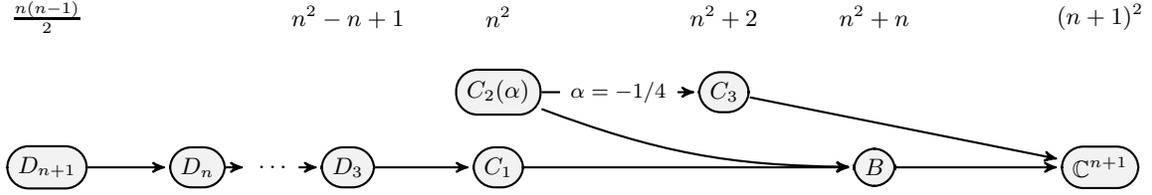
\begin{center}

\begin{tikzpicture}[->,>=stealth',shorten >=0.05cm,auto,node distance=1cm,
                    thick,main node/.style={rectangle,draw,fill=gray!10,rounded corners=1.5ex,font=\sffamily \scriptsize \bfseries },rigid node/.style={rectangle,draw,fill=black!20,rounded corners=1.5ex,font=\sffamily \scriptsize \bfseries },style={draw,font=\sffamily \scriptsize \bfseries }]

\node (0) {};

\node[main node]  (Dn1) [below   of=0]      {$D_{n+1}$};
\node (0a1) [right  of=Dn1]      {};

\node  (dn11) [above of=Dn1]      {};
\node  (dn111) [above of=dn11]      {$\frac{n(n-1)}{2}$};

\node[main node]  (Dn) [right  of=0a1]      {$D_{n}$};

\node  (dn21) [above of=Dn]      {};
%\node  (dn211) [above of=dn21]      {$n+1+\frac{(n-1)(n-2)}{2}$};

\node   (D) [right  of=Dn]      {$ \ldots$};

\node[main node]  (D3) [right  of=D]      {$D_{3}$};
\node (0d) [right  of=D3]      {};

\node  (dn31) [above of=D3]      {};
\node  (dn311) [above of=dn31]      {$n^2-n+1$};

\node[main node]  (C1) [right  of=0d]      {$C_{1}$};

\node (0aa) [right  of=C1]      {};

\node (0aaa) [right  of=0aa]      {};

\node (0aaaa) [right  of=0aaa]      {};
\node (0aaaaa) [right  of=0aaaa]      {};

\node[main node]  (B1) [right  of=0aaaaa]      {$B$};

\node  (b11) [above of=B1]      {};
\node  (b111) [above of=b11]      {$n^2+n$};

\node[main node]  (C3) [above of=0aaaa]      {$C_{3}$};
\node  (C3a) [right of=C1]      {};
\node  (c31) [below of=0aaaa]      {};
\node  (c311) [above of=C3]      {$n^2+2$};

\node[main node]  (C2) [above of=C1]      {$C_{2}(\alpha)$};

%\node[main node]  (B2) [below of=C1]      {$B_{2}$};
\node  (C2a) [right of=B1]      {};

\node  (b21) [above of=C2]      {$n^2$};

\node  (C2a2) [right of=C2a]      {};

\node[main node] (A1) [right  of=C2a2]       {$\mathbb{C}^{n+1}$};

\node  (A11) [above of=A1]      {};
\node  (A111) [above of=A11]      {$(n+1)^2$};

\path[every node/.style={font=\sffamily\small}]

% (A11a1a)  edge node[above, fill=white] {\tiny $\alpha=1,\beta=2$} (C)

(Dn1) edge  (Dn)  
(Dn)  edge  (D)  
(D)   edge  (D3)  
(D3)  edge  (C1)  
(C1)  edge  (B1)  
(B1)  edge  (A1)  

%(B2) edge [bend right=-10] (B1)  
%(C2) edge  (B1)  

(C3) edge  (A1)  
(C2) edge [bend right=10] (B1)  

(C2) edge  [bend right=0] node[above=0, right=-23, fill=white]{\tiny  $\alpha=-1/4$  } (C3)   ;        
\end{tikzpicture}

%{\bf Figure.}  The graph of primary degenerations for $(n+1)$-dimensional $n$-Lie algebras.

\end{center}

\begin{Proof}
Let us first construct all primary degenerations that are required for the assertion of the theorem.
For each degeneration $\mathbb{A} \to  \mathbb{B}$, we will give a parameterized basis $\mathcal{E}= \{E_i^t\}$, calculate the structure constants of $\mathbb{A}$ in this basis and it will be easy to see that substitution $t=0$ gives the structure constants of the algebra $\mathbb{B}$.

\begin{enumerate}
\item $C_1 \to B.$  

Parameterized basis:
$$\{E_1^t=te_1,  E_2^t=e_2, E_3^t=te_3,  E_4^t=e_4, \ldots, E_{n+1}^t=e_{n+1}\}.$$

Structure constants:
\begin{align*}
[E_2^t, E_3^t, \ldots, E_{n+1}^t]&= [e_2, te_3, e_4, \ldots, e_{n+1}]=te_1=E_1^t,\\
[E_1^t, E_3^t, \ldots, E_{n+1}^t]&= [te_1, te_3, e_4, \ldots, e_{n+1}]=t^2e_2=t^2E_2^t,\\
[E_1^t, \ldots, E_{i-1}^t,E_{i+1}^t, \ldots, E_{n+1}^t]&=[te_{1}, e_2, \ldots, e_{i-1},e_{i+1}, \ldots, e_{n+1}]=0 \mbox{ for }i>2.
\end{align*}

\item $C_2(\alpha) \to B.$ 

Parameterized basis:
$$\{E_1^t=te_2,  E_2^t=e_1, E_3^t=te_3, E_4^t=e_4, \ldots, E_{n+1}^t=e_{n+1}\}.$$

Structure constants:
\begin{align*}
[E_2^t, E_3^t, \ldots, E_{n+1}^t]&= [e_1, te_3, e_4, \ldots, e_{n+1}]=te_2=E_1^t,\\
[E_1^t, E_3^t, \ldots, E_{n+1}^t]&= [te_2, te_3, e_4, \ldots, e_{n+1}]=t^2(\alpha e_1+e_2)=tE_1+\alpha t^2 E_2^t,\\
[E_1^t, E_2^t, E_4^t,\ldots, E_{n+1}^t]&=[te_{1}, e_2,e_4, \ldots, e_{n+1}]=0,\\
[E_1^t, \ldots, E_{i-1}^t,E_{i+1}^t, \ldots, E_{n+1}^t]&=[te_{1}, e_2,te_3,e_4, \ldots, e_{i-1},e_{i+1}, \ldots, e_{n+1}]=0 \mbox{ for }i>3.
\end{align*}

\item $C_2(-1/4) \to C_3.$ 

Parameterized basis:
$$\{ E_1^t=e_1-2e_2,  E_2^t=-2te_2, E_3^t=2e_3,  E_4^t=e_4, \ldots, E_{n+1}^t=e_{n+1}\}.$$

Structure constants:
\begin{align*}
[E_2^t, E_3^t, \ldots, E_{n+1}^t]&= [-2te_2, 2e_3, e_4, \ldots, e_{n+1}]=te_1-4te_2=tE_1^t+E_2^t,\\
[E_1^t, E_3^t, \ldots, E_{n+1}^t]&= [e_1-2e_2, 2e_3, e_4, \ldots, e_{n+1}]=2e_2+e_1-4e_2=E_1^t,\\
[E_1^t, E_2^t, E_4^t,\ldots, E_{n+1}^t]&=[e_1-2e_2, -2te_2,e_4, \ldots, e_{n+1}]=0,\\
[E_1^t, \ldots, E_{i-1}^t,E_{i+1}^t, \ldots, E_{n+1}^t]&=[e_1-2e_2, -2te_2,te_3,e_4, \ldots, e_{i-1},e_{i+1}, \ldots, e_{n+1}]=0 \mbox{ for }i>3.
\end{align*}

\item $D_{3} \to C_1.$ 

Parameterized basis:
$$\{E_1^t=te_1,  E_2^t=te_2, E_{3}^t=e_3,  \ldots,  E_{n+1}^t=e_{n+1}\}.$$

Structure constants:\\
\begin{align*}
[E_2^t, E_3^t, \ldots, E_{n+1}^t]&= [te_2, e_3, \ldots, e_{n+1}]=  te_1=E_1^t,\\
[E_1^t, E_3^t, \ldots, E_{n+1}^t]&= [te_1, e_3, e_4, \ldots, e_{n+1}]=te_2=E_2^t,\\
[E_1^t, E_2^t, E_4^t, \ldots, E_{n+1}^t]&=[te_1, te_2, e_4,  \ldots, e_{n+1}]=t^2e_3=t^2E_3^t,\\
[E_1^t, \ldots, E_{i-1}^t,E_{i+1}^t, \ldots, E_{n+1}^t]&=[te_1, te_2,e_3,\ldots, e_{i-1},e_{i+1}, \ldots, e_{n+1}]=0 \mbox{ for }i>3.
\end{align*}

\item $D_{r} \to D_{r-1}, r>3.$ 

Parameterized basis:
$$\{E_1^t=te_1,  E_2^t=te_2, \ldots,  E_{r-1}^t=te_{r-1},  E_{r}^t=t^{3-r}e_{r},  E_{r+1}^t=e_{r+1}, \ldots,  E_{n+1}^t=e_{n+1}\}.$$

Structure constants:
\begin{align*}
[E_1^t, \ldots, E_{i-1}^t,E_{i+1}^t, \ldots, E_{n+1}^t]&= 
[te_1, te_2, \ldots, te_{i-1},te_{i+1}, \ldots,te_{r-1}, t^{3-r}e_{r}, e_{r+1}, \ldots, e_{n+1}]\\
&=  te_i=E_i^t\mbox{ for } 1 \leq i \leq r-1,\\
[E_1^t, E_2^t, \ldots, E_{r-1}^t,E_{r+1}^t, \ldots, E_{n+1}^t]&= 
[te_1, te_2, \ldots, te_{r-1}, e_{r+1}, \ldots, e_{n+1}]=  t^{r-1} e_r=t^{2r-4}E_r^t,\\
[E_1^t, \ldots, E_{i-1}^t,E_{i+1}^t, \ldots, E_{n+1}^t]&= 
[te_1, te_2, \ldots, te_{r-1}, t^{3-r}e_{r}, e_{r+1}, \ldots,e_{i-1},e_{i+1}, \ldots, e_{n+1}]\\
&=0\mbox{ for } r+1 \leq i \leq n+1.
\end{align*}
\end{enumerate}

Now we prove the required primary non-degenerations. It is not difficult to calculate the $(1,1)$-invariants of the algebras $C_2(\alpha)$, $C_3$ and $D_{n+1}$. One can show that
$$
c_{1,1}\big(C_2(\alpha)\big)=\frac{1}{2\alpha+1},\,c_{1,1}(C_3)=2\mbox{ and }c_{1,1}(D_{n+1})=0.
$$
Now we immediately see that $D_{n+1}\not\to C_2(\alpha),C_3$ for any $\alpha\in\mathbb{C}$ and $C_2(\alpha)\not\to C_3$ for $\alpha\not=-1/4$.

It remains to show that $C_3\not\to B$. Note that in the case of Lie algebras this assertion is valid by the automorphism group dimension argument.
Let us consider the set $\mathcal{R}$ of algebra structures $\mu$ on the space $\mathbb{V}$ with the basis $f_1,\dots,f_{n+1}$ defined by the equality
$$
\mathcal{R}=\left\{\mu\in \mathfrak{Fil}^n_{n+1}\mid
\mu(x\otimes X)\in \langle x\rangle \mbox{ for any }x\in\langle f_n,f_{n+1}\rangle\mbox{ and }X\in\mathbb{V}^{\otimes(n-1)}
\right\},
$$
where $\langle y_1,\dots,y_l\rangle$ denotes the subspace of $\mathbb{V}$ generated by $y_1,\dots,y_l\in\mathbb{V}$. It is not difficult to see that $\mathcal{R}$ is a closed subset of $\mathfrak{Fil}^n_{n+1}$ that is stable under the action of the subgroup of $GL(\mathbb{V})$ formed by automorphisms having lower triangular matrices in the basis $f_1,\dots,f_{n+1}$. To show that $C_3$ can be represented by a structure from $\mathcal{R}$, it is enough to take $f_i=e_{n+2-i}$ for $1\le i\le n+1$. On the other hand, if $B(x\otimes X)\in \langle x\rangle$ for any $X\in\mathbb{V}^{\otimes(n-1)}$, then $x\in\langle e_1\rangle$. But for any $\mu\in\mathcal{R}$ the set of such elements $x$ that $\mu(x\otimes X)\in \langle x\rangle$ for any $X\in\mathbb{V}^{\otimes(n-1)}$ contains a subspace of dimension not less than two. Thus, $C_3\not\to B$ by Lemma \ref{main}.
\end{Proof}

Note that, for $\alpha\not=0$, the algebra $C_2(\alpha)$ is isomorphic to the algebra whose multiplication table is determined by the nonzero products $[e_2, \ldots,  e_{n+1}]=e_1+\frac{1}{\alpha}e_2$ and $[e_1, e_3,  \ldots, e_{n+1}]=e_2$. To see this, it is sufficient to replace $e_2$ and $e_3$ by $\frac{e_2}{\alpha}$ and $\frac{e_3}{\alpha}$ respectively in the standard basis. Then it is clear that $C_1$ belongs to $\overline{\{O\big(C_2(\alpha)\big)\}_{\alpha\in\mathbb{C}}}$. Thus, we get the description of irreducible components of $\mathfrak{Fil}^n_{n+1}$ presented in the next corollary.

\begin{corollary} The irreducible components  of   $\mathfrak{Fil}^n_{n+1}$ are
\begin{align*}
\mathcal{C}_1&=\overline{\{O\big(C_2(\alpha)\big) \}_{\alpha \in\mathbb{C}}} = \{\mathbb{C}^{n+1},  B, C_1, C_2(\alpha), C_3 \},\\
\mathcal{C}_2&=\overline{O(D_{n+1})} =  \{\mathbb{C}^{n+1}, B,C_1, D_3, \ldots, D_{n+1}  \}.
\end{align*}
In particular, $D_{n+1}$ is the unique rigid algebra in the variety $\mathfrak{Fil}^n_{n+1}$.
\end{corollary}

Let us recall that a $k$-dimensional algebra $\mathbb{A}$ has level $m$ if there is a chain
$$\mathbb{A}=\mathbb{A}_0\to \mathbb{A}_1\to\cdots\to \mathbb{A}_m=\mathbb{C}^k$$
 of primary degenerations of length $m$ and there is no such a chain of length $m+1$. In particular, a nonzero structure $\mu\in  Hom(\mathbb{V}^{\otimes n},\mathbb{V})$ with $k$-dimensional space $\mathbb{V}$ has level one if and only if $\overline{O(\mu)}=\{\mu,\mathbb{C}^k\}$.

\begin{corollary}\label{lone}
Algebras $B$ and $C_3$ have level one: 
$\overline{O(B)} = \{ B, \mathbb{C}^{n+1} \}$ and $\overline{O(C_3)} = \{ C_3 , \mathbb{C}^{n+1} \}.$
\end{corollary}

\section{Open problems}
The binary algebras of the first and the second level were classified in \cite{khud13,kv18}.
The author of \cite{gorb93} defined also the notion of an infinite level and considered the problem of classification of anticommutative algebras of small infinite levels. This notion is much easier in the sense that the infinite level of an algebra can be easily expressed in terms of the usual level.

In Corollary \ref{lone} we have showed that there exist two $(n+1)$-dimensional $n$-ary Filippov algebras of level one. 
It would be interesting to know the answers to more general questions on the levels of $n$-ary algebras.

\

{\bf Problem.}
Classify $n$-ary algebras of the first usual and infinite levels.

\

It would be interesting also to solve this problem for some varieties of $n$-ary algebras. For example, for (anti)commutative $n$-ary algebras or even for Filippov algebras of an arbitrary dimension.

%\newpage

\end{document}